\documentclass[english]{article}
\usepackage[T1]{fontenc}
\usepackage[latin9]{inputenc}
\usepackage{float}
\usepackage{amsmath}
\usepackage{amssymb}
\usepackage{stackrel}
\usepackage{esint}
\usepackage{nameref}

\makeatletter

\floatstyle{ruled}
\newfloat{algorithm}{tbp}{loa}
\providecommand{\algorithmname}{Algorithm}
\floatname{algorithm}{\protect\algorithmname}

\usepackage{graphicx}

\usepackage{xr}
\usepackage{hyperref}
\usepackage[capitalise]{cleveref}

\usepackage{tikz}
\usepackage{pgfplots}
\usepackage{pgfplotstable}
\usepgfplotslibrary{units}
\usepackage{booktabs}
\usepackage{authblk}
\usepackage[english]{babel}

\makeatother

\usepackage{babel}
\begin{document}

\title{An Elastic Energy Minimization Framework for Mean Contour Calculation}

\author[$1$]{Jozsef Molnar} 
\author[$2$]{Michael Barbier}
\author[$2$]{Winnok H. De Vos}
\author[$1,3$]{Peter Horvath}

\affil[$1$]{Synthetic and Systems Biology Unit\\
Biological Research Centre\\
Hungarian Academy of Sciences\\
Szeged, Hungary\\

e-mail: jmolnar64@digikabel.hu, horvath.peter@brc.mta.hu}

\affil[$2$]{Laboratory of Cell Biology and Histology\\
Department of Veterinary Sciences\\
Antwerp University\\
Universiteitsplein 1,2610 Antwerpen, Belgium

e-mail: (Michael.Barbier, Winnok.DeVos)@uantwerpen.be}

\affil[$3$]{Institute for Molecular Medicine Finland (FIMM)\\
University of Helsinki\\
Helsinki, Finland}

\maketitle
\begin{abstract}
In this paper we propose a contour mean calculation and interpolation
method designed for averaging manual delineations of objects performed
by experts and interpolate 3D layer stack images. The proposed method
retains all visible information of the input contour set: the relative
positions, orientations and size, but allows invisible quantities
- parameterization and the centroid - to be changed. The chosen representation
space - the position vector rescaled by square root velocity - is
a real valued vector space on which the imposed $\mathbb{L}^{2}$
metric is used to define the distance function. With respect to this
representation the re-parameterization group acts by isometries and
the distance has well defined meaning: the sum of the central second
moments of the coordinate functions. To identify the optimal re-parameterization
system and proper centroid we use double energy minimization realized
in a variational framework.
\end{abstract}

\section{Inroduction}

A specifically designed mathematical framework for two practical problems:
contour averaging and interpolation is proposed and examined in this
paper. 

Object delineation is an important annotation step to create training
data set for the supervised machine learning methods designed for
object segmentation. Histopathology images, however rarely provide
definite unambiguous object boundaries, often the delineations performed
by experts do not agree. One plausible approach to create meaningful
annotation samples is to accept the mean of many recommendations excluding
some outliers. This approach requires well defined, meaningful metrics
on the space of contours.

The resolution of several microscopy techniques in the direction of
focusing (direction $Z$) is usually a magnitude less than the resolution
of the stack images. Interpolation needs to be carried out in a principled
manner to achieve good estimation for the accurate 3D measurements
of the object physical quantities, such as surface area or volume.
Interpolation can also be useful tool to track the progression of
lesions in various diagnostic images.

The proposed method is designed to keep all visible information encoded
in the set of the constituent contours including their relative displacement,
hence essentially position vector based. The description of the contours
by preselected position vector set (landmark points) is the approach
of the early shape analysis techniques (with the identification of
shape manifolds of $k$-points and the imposed Riemannian metric see
for example \cite{Kendall84shapemanifolds}). On the other hand, the
predetermined sampling strategy of the landmark points is related
to the fixed parameterization of the contours. In the proposed model
this restriction is relaxed and some tools borrowed from the elastic
shape analysis \cite{Mio2007} are used. The chosen contour representation
is the position vector rescaled by square root velocity that - wrt
a properly defined centroid - provides covariant description whilst
retain all contextual information. It can be considered as the combination
of the landmark based and the Square Root Velocity Function (SRVF)
\cite{joshi2007novel}\cite{10.1007/978-3-540-74198-5_30}\cite{srivastava2011shape}
representations (for which the analysis of the existence of the optimal
reparameterization is found in \cite{doi:10.1137/15M1014693}). The
proposed representation and the associated $\mathbb{L}^{2}$ metric
are exhaustively examined in this paper mentioning some perspective
generalizations. References to the SRVF are also provided wherever
informative/relevant.

The structure of the paper is the following. Section \ref{sec:The-model}
presents the framework including numerical methods. Section \ref{sec:Illustrative-examples}
is dedicated to illustrative interpolation examples, section \ref{sec:Conclusion}
concludes the paper with discussion and outlook. Appendices contain
important proofs and derivations.

\section{The contour averaging framework\label{sec:The-model}}

We consider simple, planar contours used to delineate objects to be
closed, continuous, one-parameter ($t\in\left[0,T\right]$) family
objects with winding number one. From now on we simple refer them
as 'contours'. The principal representations of contours are often
given by position vector wrt some standard basis $\mathbf{i},\,\mathbf{j}$
as $\mathbf{r}\left(t\right)=x\left(t\right)\mathbf{i}+y\left(t\right)\mathbf{j}$,
$\mathbf{r}\left(0\right)=\mathbf{r}\left(T\right)$ where $x\left(t\right),\, y\left(t\right)$
are the coordinate functions. The set of contours used to calculate
their mean is referred as contour system.

To develop a framework for efficient contour mean calculation, first
we assess some natural conditions to be fulfilled by any model developed
for this purpose:
\begin{description}
\item [{A)}] Keep all visible information (relative positions, rotations,
size) of the constituents, optimize only for non-visible ones
\item [{B)}] The mean contour derived from the system needs to be invariant
for its constituents common translation, rotation, scaling (\textit{i.e.}
the mean of the transformed system is transformed in the same manner
as the constituents)
\item [{C)}] The result of the mean determination must be independent of
the parameterization of the constituents
\end{description}
The position vector representation obviously satisfies condition \textbf{A}.
It also satisfies condition \textbf{B}, if the basis $\mathbf{i},\,\mathbf{j}$
is determined by the system itself. Condition \textbf{C} however cannot
be fulfilled by this representation. One of the possibility to get
simple parameterization-invariant representation - known from the
shape analysis literature - is the choice of the square root velocity
function (SRVF) \cite{srivastava2011shape}. SRVF however, does not
retain the relative translation information. For this reason we use
the combination of the position vector and the SRVF: the position,
rescaled by square root velocity (Rescaled Position by Square Velocity
or RPSV):
\begin{equation}
\mathbf{q}\left(t\right)\doteq\mathbf{r}\left(t\right)\sqrt{\left|\dot{\mathbf{r}}\left(t\right)\right|}\,.
\end{equation}
The points of the position vector $\mathbf{r}\left(t\right)$ and
its RPSV representation $\mathbf{q}\left(t\right)$ lie in the same
direction $\mathbf{u}\left(t\right)=\frac{\mathbf{r}\left(t\right)}{\left|\mathbf{r}\left(t\right)\right|}=\frac{\mathbf{q}\left(t\right)}{\left|\mathbf{q}\left(t\right)\right|}$,
hence reproducing the contour (its position vector) requires the determination
of its length $\left|\mathbf{r}\left(t\right)\right|$ at each parameter
value $t$. This can be done iteratively using the Newton\textendash Raphson
method (see \nameref{sec:Appendix-C}).

\subsection{Properties of the representation\label{sub:Representetion-properties}}

Position vector representation $\mathbf{r}$: $\left[0,T\right]\rightarrow\mathbb{R}^{2}$
is the vector space of coordinate function duplets, so its reparameterization
$\mathbf{q}$.%
\footnote{Representations $\mathbf{r}$ and$\mathbf{q}$ are considered as as
two parameterization of the underlying space of function duplets.%
} Equipped with the inner product
\begin{equation}
\left\langle \mathbf{q}_{1},\mathbf{q}_{2}\right\rangle \doteq\oint\mathbf{q}_{1}\left(t\right)\cdot\mathbf{q}_{2}\left(t\right)dt
\end{equation}
(where $\mathbf{q}_{1}\left(t\right)\cdot\mathbf{q}_{2}\left(t\right)$is
the dot product of the position vectors given at parameter value $t$)
and the distance function based on the $\mathbb{L}^{2}$ norm $\left\Vert \mathbf{q}\right\Vert ^{2}=\left\langle \mathbf{q},\mathbf{q}\right\rangle $:
\begin{equation}
d^{2}\left(\mathbf{q}_{1},\mathbf{q}_{2}\right)=\oint\left(\mathbf{q}_{1}\left(t\right)-\mathbf{q}_{2}\left(t\right)\right)^{2}dt\label{eq:squared-distance}
\end{equation}
the space of the representations $\mathbf{q}$ becomes Hilbert space,
denoted by $\mathcal{H}_{q}$. With reference to the \nameref{sec:Appendix-A}
here we asses the important properties of the chosen representataion:
\begin{enumerate}
\item The squared norm $\left\Vert \mathbf{q}\right\Vert ^{2}$ of any point
in the representation space expresses the sum of the second central
moments of the (coordinate functions of) contour $\mathbf{r}\left(t\right)$,
consequently:
\item The distance function is invariant wrt the common reparameterization
of points $\mathbf{q}_{1},\,\mathbf{q}_{2}\rightarrow\mathbf{q}_{1}\circ\gamma,\,\mathbf{q}_{2}\circ\gamma$
and
\item The reparameterization group $\Gamma=\left\{ \left.\gamma\right|\gamma\left(t\right)>0\right\} $
acts by isometries wrt the chosen metric and composition, \textit{i.e.}
$d^{2}\left(\mathbf{q}_{1},\mathbf{q}_{2}\right)=d^{2}\left(\mathbf{q}_{1}\circ\gamma,\mathbf{q}_{2}\circ\gamma\right)$
\end{enumerate}
The properties above allows us to construct the mean contour in the
quotient space $\mathcal{H}_{q}/\Gamma$ satisfying the requirements
\textbf{A}, \textbf{B}, \textbf{C} stipulated at the beginning of
this section. The mean representation of a system of $n$ representations
$\mathbf{q}_{1},\,\mathbf{q}_{2},\cdots\mathbf{q}_{n}$ is defined
to be:
\begin{equation}
\mathbf{q}\left(t\right)\doteq\frac{1}{n}\stackrel[i=1]{n}{\sum}\mathbf{q}_{i}\left(t_{i}\right),\, t,t_{i}\in\left[0,T\right]\label{eq:mean-by-representation}
\end{equation}
where parameterizations $t_{i}=\gamma_{i}\left(t\right)$ are the
carefully selected points from the orbit of $\Gamma$ the for which
$\stackrel[i=1]{n}{\sum}d^{2}\left(\mathbf{q}_{i}\left(t_{i}\right)-\mathbf{q}\left(t\right)\right)$
is minimal. Formula (\ref{eq:mean-by-representation}) can be written
directly in position vector 'coordinates' of the space $\mathcal{H}_{q}$
and takes the form:
\begin{equation}
\mathbf{r}\left(t\right)\sqrt{\left|\dot{\mathbf{r}}\left(t\right)\right|}\doteq\frac{1}{n}\stackrel[i=1]{n}{\sum}\mathbf{r}_{i}\left(t_{i}\right)\sqrt{\left|\dot{\mathbf{r}}_{i}\left(t_{i}\right)\right|}\,.\label{eq:mean-by-position}
\end{equation}
Later in the paper we use mainly the direct (position vector) coordinates,
not forgetting the underlying RPSV representation. We conclude this
subsection with the statement: the mean contour is identified as the
position vector associated with the mean of the RPSV representetions.

\subsection{Mean contour as a minimization problem}

Properties 1-3 of RPSV (\nameref{sub:Representetion-properties})
enable to construct the optimal parameterization of the system of
contours in a simple way, \textit{i.e.} choosing one of the constituent
contour as 'reference contour' and calculate the optimal parameterization
of the other contours wrt it, see Lemma A1 in \nameref{sec:Appendix-A}. 

The minimization problem wrt a fixed origin (will be relaxed later),
using direct position vector coordinates can be formulated as:
\begin{equation}
\underset{\gamma_{i}}{\min}\stackrel[i=1]{n}{\sum}\oint\left(\mathbf{r}\left(t\right)\sqrt{\left|\dot{\mathbf{r}}\left(t\right)\right|}-\mathbf{r}_{i}\left(t_{i}\right)\sqrt{\left|\dot{\mathbf{r}}_{i}\left(t_{i}\right)\right|}\right)^{2}dt\label{eq:minproblem-raw}
\end{equation}
where $\mathbf{r}\left(t\right)\sqrt{\left|\dot{\mathbf{r}}\left(t\right)\right|}$
stands for the mean contour (\ref{eq:mean-by-position}), $t_{i}=\gamma_{i}\left(t\right)$.
As analysed in \nameref{sec:Appendix-A}, the solution (system of
$\gamma_{i}$) that provides the minimum distances between the constituents
$\mathbf{r}_{i}\left(t_{i}\right)\sqrt{\left|\dot{\mathbf{r}}_{i}\left(t_{i}\right)\right|}$
can be determined pairwise wrt a reference contour (say $\mathbf{r}_{1}$
without loss of generality)
\begin{equation}
\underset{\gamma_{k}}{\min}\oint\left(\mathbf{r}_{1}\left(t\right)\sqrt{\left|\dot{\mathbf{r}}_{1}\left(t\right)\right|}-\mathbf{r}_{k}\left(t_{k}\right)\sqrt{\left|\dot{\mathbf{r}}_{k}\left(t_{k}\right)\right|}\right)^{2}dt\label{eq:minproblem-pairwise}
\end{equation}
 as the solution of the Euler-Lagrange equations assiciated with them:
\begin{equation}
\dot{\mathbf{r}}_{1}\cdot\mathbf{r}_{k}-\dot{\mathbf{r}}_{k}\cdot\mathbf{r}_{1}+\frac{1}{2}\left(\varGamma_{k}-\varGamma_{1}\right)=0\label{eq:Euler-Lagrange-raw}
\end{equation}
where the dot over the position vectors stands for the derivatives
wrt the parameter $t$, \textit{i.e.} $\dot{\mathbf{r}}_{k}\equiv\frac{d\mathbf{r}_{k}}{dt}=\dot{\gamma}_{k}\frac{d\mathbf{r}_{k}}{d\gamma_{k}}$,
$k=2,\ldots n$, (note: since $\mathbf{r}_{1}$ is chosen as reference
contour $\mathbf{r}_{1}\left(\gamma_{1}\left(t\right)\right)\equiv\mathbf{r}_{1}\left(t\right)$
in (\ref{eq:minproblem-raw}) ) and 
\begin{equation}
\varGamma_{i}=\frac{\dot{\mathbf{r}}_{i}\cdot\ddot{\mathbf{r}}_{i}}{\left|\dot{\mathbf{r}}_{i}\right|^{2}},\, i=1,\ldots n,\label{eq:Christoffel-divergences}
\end{equation}
are the 'Christoffel divergences' of the parameterization. As expected,
the solution for the minimization problem (\ref{eq:minproblem-raw})
is given by the system $\gamma_{k}$ determined pairwise, using Euler-Lagrange
equation (\ref{eq:Euler-Lagrange-raw}), see also \nameref{sec:Appendix-B}.
Notes:
\begin{enumerate}
\item Euler-Lagrange equation (\ref{eq:Euler-Lagrange-raw}) retains its
form wrt any basis, albeit the resulting system of the optimal parameterization
is dependent on the chosen basis; we will address this problem in
section \nameref{sub:Proper-centroid}
\item (apart from a proportionality factor) Euler-Lagrange equation (\ref{eq:Euler-Lagrange-raw})
does not depend explicitly on the reparameterization function $\gamma_{k}$
\item Christoffel divergences $\varGamma=\frac{d\ln\left|\dot{\mathbf{r}}\right|}{dt}$
can be interpreted as the change of 'elastic stretching' along the
contours; indeed the quantity $\ln\left|\dot{\mathbf{r}}\right|$
has prominent role in definition of elastic shape metrics in \cite{Mio2007}
\item assuming $\mathbf{r}_{1}$ is uniformly parameterized in arc length:$\varGamma_{1}=0$
and $\varGamma_{k}=2\left(\dot{\mathbf{r}}_{k}\cdot\mathbf{r}_{1}-\dot{\mathbf{r}}_{1}\cdot\mathbf{r}_{k}\right)$,
$k=2,\ldots n$ determine the elastic stretching/compression
\item from a different point of view, Eq.(\ref{eq:minproblem-pairwise})
can be considered as 'dissimilarity measure' between contours 
\item as expected, exactly same Euler-Lagrange equations (\ref{eq:Euler-Lagrange-raw})
are associated with the similarity maximization $\underset{\gamma_{k}}{\max}\oint\mathbf{r}_{1}\left(t\right)\cdot\mathbf{r}_{k}\left(t_{k}\right)\sqrt{\left|\dot{\mathbf{r}}_{1}\left(t\right)\right|\left|\dot{\mathbf{r}}_{k}\left(t_{k}\right)\right|}dt$,
$k=2,\ldots n$ problems
\item note that in the SRVF case, the optimal reparameterization problen
can also be formulated as variational problem and its associated Euler-Lagrange
equation can be arranged to $\ddot{\mathbf{r}}_{1}\cdot\mathbf{\dot{r}}_{k}-\ddot{\mathbf{r}}_{k}\cdot\dot{\mathbf{r}}_{1}+\frac{1}{2}\left(\varGamma_{k}-\varGamma_{1}\right)=0$,
having formal similariy to equation (\ref{eq:Euler-Lagrange-raw})
with higher-order derivations applied to the first two terms.
\end{enumerate}

\subsubsection{Proper centroid\label{sub:Proper-centroid}}

To elaborate a covariant model, the origin of the standard basis $\mathbf{i},\,\mathbf{j}$
wrt the position vectors are expressed must be defined by the contour
system itself. Otherwise the mean contour would not be invariant to
the common translation of its constituents, violating requirement
\textbf{B} stated at the beginning of this section \ref{sec:The-model}.
Now assume, we have our contour system wrt some ad hoc basis and denote
the position vector wrt that basis with $\mathbf{R}_{i}\left(t\right)$,
$i=1,\ldots n$. First plausible candidate for the origin would be
the usual centroid of the system that minimizes:
\begin{equation}
\underset{\mathbf{R}_{0}}{\min}\stackrel[i=1]{n}{\sum}\oint\left(\mathbf{R}_{0}-\mathbf{R}_{i}\left(t\right)\right)^{2}\left|\dot{\mathbf{R}}_{i}\right|dt\,.\label{eq:homogen-centroid}
\end{equation}
This candidate provides covariant description, also independent of
the parameterization of the constituents (since $\left|\dot{\mathbf{R}}_{i}\right|dt=ds$,
the integration is by arc length). From now on we will refer to it
as the 'homogeneous' centroid. Problem (\ref{eq:homogen-centroid})
can be interpreted as simple extreme value problem wrt the centroid
coordinates $\mathbf{R}_{0}$ and such the condition: 
\begin{equation}
\frac{d\left(\stackrel[i=1]{n}{\sum}\oint\left(\mathbf{R}_{0}-\mathbf{R}_{i}\left(t\right)\right)^{2}\left|\dot{\mathbf{R}}_{i}\right|dt\right)}{d\mathbf{R}_{0}}\doteq0
\end{equation}
provides the following solution:
\begin{equation}
\mathbf{R}_{0}=\frac{\stackrel[i=1]{n}{\sum}\oint\mathbf{R}_{i}ds}{\stackrel[i=1]{n}{\sum}L_{i}}\label{eq:homogeneous-centroid}
\end{equation}
where $L_{i}$ stands for the length of the $i$-th contour. Adopting
the standard basis to be this homogeneous centroid, the position vectors
of the contour system wrt this basis would become $\mathbf{r}_{i}\left(t\right)=\mathbf{R}_{i}\left(t\right)-\mathbf{R}_{0}$
$i=1,\ldots n$.

However the question arises naturally: is the choice of the homogeneous
centroid 'compatible' with the minimization problem (\ref{eq:minproblem-raw})?
To decide this question, let's assume, we displace the basis from
the homogeneous centroid position with a vector $\delta\mathbf{d}$.
The position vectors are then transformed to $\mathbf{r}_{i}\left(t\right)\rightarrow\mathbf{r}_{i}\left(t\right)-\mathbf{d}$.
Now check, whether the double minimization problem, generalized from
(\ref{eq:minproblem-raw}):
\begin{equation}
E\left(\gamma_{i},\delta\mathbf{d}\right)=\underset{\gamma_{i},\delta\mathbf{d}}{\min}\stackrel[i=1]{n}{\sum}\oint\left[\left(\mathbf{r}\left(t\right)-\delta\mathbf{d}\right)\sqrt{\left|\dot{\mathbf{r}}\left(t\right)\right|}-\left(\mathbf{r}_{i}\left(t_{i}\right)-\delta\mathbf{d}\right)\sqrt{\left|\dot{\mathbf{r}}_{i}\left(t_{i}\right)\right|}\right]^{2}dt\label{eq:minproblem-fine}
\end{equation}
takes its minimum at $\delta\mathbf{d}=\mathbf{0}$. From the condition
$\frac{\partial E}{\partial\delta\mathbf{d}}=0$, one can derive:
\begin{equation}
\delta\mathbf{d}=\frac{\stackrel[i=1]{n}{\sum}S_{i}+nS-\left[\stackrel[i=1]{n}{\sum}\oint\left(\mathbf{r}\left(t\right)+\mathbf{r}_{i}\left(t_{i}\right)\right)\sqrt{\left|\dot{\mathbf{r}}\left(t\right)\right|\left|\dot{\mathbf{r}}_{i}\left(t_{i}\right)\right|}dt\right]}{\stackrel[i=1]{n}{\sum}L_{i}+nL-\left[2\stackrel[i=1]{n}{\sum}\oint\sqrt{\left|\dot{\mathbf{r}}\left(t\right)\right|\left|\dot{\mathbf{r}}_{i}\left(t_{i}\right)\right|}dt\right]}\label{eq:proper-centroid}
\end{equation}
where notations $S_{i}=\oint\mathbf{r}_{i}\left(t_{i}\right)\sqrt{\left|\dot{\mathbf{r}}_{i}\left(t_{i}\right)\right|}dt$,
$S=\oint\mathbf{r}\left(t\right)\sqrt{\left|\dot{\mathbf{r}}\left(t\right)\right|}dt$
and the lengths of the constituents $L_{i}=\oint\left|\dot{\mathbf{r}}_{i}\left(t_{i}\right)\right|dt$
and the mean contour $L=\oint\left|\dot{\mathbf{r}}\left(t\right)\right|dt$
are introduced. Now one can notice that in general, the optimal displacement
of the homogeneous centroid wrt minimization problem (\ref{eq:minproblem-fine})
is not zero vector due to the parameterization dependent terms emphasized
in brackets in (\ref{eq:proper-centroid}). This issue obviously stem
from the fact that the optimally parameterized contour system consists
of non-uniformly parameterized (in arc length sense) 'inhomogeneous'
contours. From now on we refer the centroid that satisfies the double
minimization problem (\ref{eq:minproblem-fine}) as proper centroid.

The optimal centroid and parameterization system are interdependent:
wrt a fixed basis a unique optimal reparameterization system can be
calculated which in turn determines the location of the proper centroid;
on the other hand in general (unless $\delta\mathbf{d}=\mathbf{0}$
by (\ref{eq:proper-centroid})) the optimal parameterization system
is dependent on the choice of the standard basis. This interdependency
leads to an iterative solution which is discussed in details in the
next section. The optimal reparameterization system and the proper
centroid are determined alternately. Using this approach, equation
(\ref{eq:proper-centroid}) can be simplified as follows. In the first
step the optimal reparameterization system is determined wrt the momentary
centroid, then the mean contour is calculated (\ref{eq:mean-by-position})
and reconstructed. Substituting the mean $\mathbf{r}\sqrt{\left|\dot{\mathbf{r}}\right|}\doteq\frac{1}{n}\stackrel[i=1]{n}{\sum}\mathbf{r}_{i}\sqrt{\left|\dot{\mathbf{r}}_{i}\right|}$
to $nS$ in (\ref{eq:proper-centroid}), two terms are eliminated
from the enumerator. After some rearrangement (both the enumerator
and the denominator) we arrive to a simple expression: 
\begin{equation}
\delta\mathbf{d}=\frac{\stackrel[i=1]{n}{\sum}\oint\left(\mathbf{r}_{i}\sqrt{\left|\dot{\mathbf{r}}_{i}\right|}-\mathbf{r}\sqrt{\left|\dot{\mathbf{r}}\right|}\right)\sqrt{\left|\dot{\mathbf{r}}_{i}\right|}dt}{2\stackrel[i=1]{n}{\sum}\oint\frac{\left|\dot{\mathbf{r}}_{i}\right|+\left|\dot{\mathbf{r}}\right|}{2}-\sqrt{\left|\dot{\mathbf{r}}_{i}\right|\left|\dot{\mathbf{r}}\right|}dt}\,.\label{eq:proper-centroid-1}
\end{equation}
In the denominator, the integrand is the sum of the differences of
the arithmetic and geometric means of the corresponding elementary
arc lengths $\frac{ds+ds_{i}}{2}$ and $\sqrt{ds_{i}ds}$ respectively
(using the $ds=\left|\dot{\mathbf{r}}\right|dt$ identity). The denominator
therefore can be zero only if the lengts of all the corresponding
elementary arc segments are identical, the case possible only if the
constituent contours are all identical.

\subsection{Numerical methods\label{sub:Numerical-methods}}

As in the case of shape analysis, the calculation of the mean contour
requires iterative solutions: a double iteration for determination
of the optimal reparameterization system and the proper centroid defined
by (\ref{eq:minproblem-fine}), then one for the reconstruction of
the contour from its RPSV representetion. The components are the following. 
\begin{enumerate}
\item Reparameterization of the system
\item Mean calculation
\item Reconstruction of the mean contour from its representation
\item Proper centroid calculation
\end{enumerate}

\subsubsection*{Reparamaterization\label{sub:Reparamaterization}}

The identification of the optimal reparameterization system (\ref{eq:minproblem-fine})
requires the calculation of $n-1$ pairwise reparameterization wrt
a reference contour. The gradient descent equations are 
\begin{eqnarray}
\frac{\partial\gamma_{k}}{\partial\tau} & = & -\dot{\mathbf{r}}_{1}\cdot\mathbf{r}_{k}+\mathbf{r}_{1}\cdot\dot{\mathbf{r}}_{k}-\frac{1}{2}\mathbf{r}_{1}\cdot\mathbf{r}_{k}\left(\varGamma_{1}-\varGamma_{k}\right),\label{eq:Euler-Lagrange-fine}
\end{eqnarray}
where $\tau$ is the 'artifical' time and the Christoffer divergences
are defined by (\ref{eq:Christoffel-divergences}). These equations
are to be solved in the contour space. Two methodologies are possible
to determine the optimal parameterization. In the first (recommended)
case, after each iteration, the points are redistributed moving them
to their new physical position determined by $\delta\gamma_{k}^{\left(i\right)}$
($i=1,...N$ is the iteration index) along the (static) contours $\mathbf{r}_{k}$.
Derivatives $\dot{\mathbf{r}}_{k}$ are calculated from the momentary
positions of the contour points. Note that in the discrete approximation
of contours, uniform distribution wrt parameter value $t$ can be
assumed without loss of generality (that is the parameter values assigned
to the neighboring points differ from each-other with same $\varDelta t$
everywhere). In this case in the parameter space the parameter values
associated with the (moving) points remain constant, albeit their
arclength parameters change in general. The final diffeomorphism $\gamma_{k}$
is then the composition of the sequence of consequtive approximate
diffeomorphisms $\delta\gamma_{k}^{\left(i\right)}$, that is (assuming
overall $N$ iterations) $\gamma_{k}=\delta\gamma_{k}^{\left(N\right)}\circ\cdots\circ\delta\gamma_{k}^{\left(2\right)}\circ\delta\gamma_{k}^{\left(1\right)}$.
See Algorithm \ref{alg:Compute-pairwise-optimal}. In the second case
$\gamma_{k}$ is updated after all iterations with the points physical
position retained at their initial position. This approach however,
requires the calculations of the derivatives wrt the momentary $\gamma_{k}^{\left(i\right)}$
using explicite formulae for the derivatives (\textit{i.e.} $\frac{d}{dt}=\dot{\gamma}\frac{d}{d\gamma}$).
The first methology has the advantages a) at each iteration step $\delta\gamma_{k}^{\left(i\right)}$
needs to be determined wrt the identity diffeomorphism $\gamma\left(t\right)\equiv t$
b) usually, there is no real need for the explicit determination of
the final diffeomorphism, only the final point distribution we end
up with the first methology and c) it can be efficiently implemented
using a high resolution lookup table for the positions along the contours.

\begin{algorithm}
\protect\caption{Compute pairwise optimal reparameterization\label{alg:Compute-pairwise-optimal}}

\begin{enumerate}
\item Initialize the position vectors $\mathbf{r}_{k}$, $k=1\ldots n$
of the contour set wrt the homogeneous centroid (\ref{eq:homogeneous-centroid}).
Establish the initial discrete point set along the contours with same
number of points (can be uniformly distributed in arc length); Set
the iteration counter $i=1$; Set $\delta\gamma_{k}^{\left(0\right)}=t$,
$k=2\ldots n$ (\textit{i.e.} $\mathbf{r}_{1}$ is selected as reference)
\item Calculate one step towards ($\delta\gamma_{k}^{\left(i\right)}$,
$k=2\ldots n$) the optimal point distribution system using gradient
descent equations (\ref{eq:Euler-Lagrange-fine})
\item Update the points along contours $\mathbf{r}_{k}$, $k=2\ldots n$,
using the calculated valuea $\delta\gamma_{k}^{\left(i\right)}$,
$k=2\ldots n$
\item Update the diffeomorphism set $\gamma_{k}^{\left(i\right)}=\delta\gamma_{k}^{\left(i\right)}\circ\gamma_{k}^{\left(i-1\right)}$
\item Exit if all $\delta\gamma_{k}^{\left(i\right)}$ (wrt its $\gamma_{k}^{\left(i\right)}$)
is small; Otherwise set $\delta\gamma_{k}^{\left(i\right)}=t$, set
$i=i+1$ and repeat from 2
\end{enumerate}
\end{algorithm}

\subsubsection*{Mean calculation\label{sub:Mean-calculation}}

Given the optimal reparameterization system, the mean is calculated
using the closed form equation (\ref{eq:mean-by-representation}).

\subsubsection*{Recontruction\label{sub:Recontruction}}

Reconstruction is made by the Newton\textendash Raphson method, solving
a sparse linear equation system in each iteration $\mathbf{A}^{\left(i\right)}\mathbf{x}^{\left(i+1\right)}=\mathbf{b}^{\left(i\right)}$
($i$ is the iteration index) with coefficient matrix, ray length
approximation of the position vector and constant vector all defined
in \nameref{sec:Appendix-C} by formulae (\ref{eq:coeffitient-matrix}),
(\ref{eq:ray-length}), (\ref{eq:constant-vector}) respectively.

\subsubsection*{Proper centroid\label{sub:Proper-centroid-1}}

Proper centroid for the momentary parameterization system is calculated
using the closed form formula (\ref{eq:proper-centroid-1}). Once
the (better) displacement $\delta\mathbf{d}^{\left(j\right)}$ is
determined all constituent contours have to be updated such as $\mathbf{r}_{k}\rightarrow\mathbf{r}_{k}+\delta\mathbf{d}^{\left(j\right)}$,
$k=1,\ldots n$ then all previous steps are to be repeated until the
minimum of the double minimization problem (\ref{eq:minproblem-fine})
is reached. The cumulative displacement of the initial (homogeneous)
centroid after $M$ iterations is the sum of the preceding (momentary)
displacements: $\stackrel[j=1]{M}{\sum}\delta\mathbf{d}^{\left(j\right)}$.

\subsubsection*{The algorithm}

Albeit the determination of the optimal reparameterization system
and the proper centroid calculation could be incorporated into one
iterative method, but the need for the mean contour calculation in
(\ref{eq:proper-centroid}) after each gradient descent step of (\ref{eq:Euler-Lagrange-fine})
would lead to sluggish computing. Therefore a double iteration procedure
is recommended: an inner (nested) loop for the optimal reparameterization
system under the assumption of centroid constancy, followed by the
centroid position updating in the outer (main) loop. 

The complete algorithm consists of the steps described above and summarized
in Algorithms \ref{alg:Compute-pairwise-optimal} (nested loop) and
\ref{alg:Solve-the-double} (main loop).

\begin{algorithm}
\protect\caption{Solve the double optimization algorithm\label{alg:Solve-the-double}}

\begin{enumerate}
\item Initialize the position vectors $\mathbf{r}_{k}$, $k=1\ldots n$
of the contour set wrt the homogeneous centroid (\ref{eq:homogeneous-centroid}).
Establish the initial discrete point set along the contours with same
number of points (can be uniformly distributed in arc length); Set
the iteration counter $j=1$; Set $\delta\mathbf{d}^{\left(j\right)}=0$.
\item In internal loop compute the optimal redistribution system of points
pairwise wrt an arbirtarily designated reference contour using gradient
descent equation (\ref{eq:Euler-Lagrange-fine}) or alternatively
compute the optimal reparameterization system $\gamma_{i}$: \nameref{sub:Reparamaterization};
see also Algorithm \ref{alg:Compute-pairwise-optimal}
\item Calculate the mean contour in the representation space RPSV: \nameref{sub:Mean-calculation}
\item Reconstruct the mean in contour space: \nameref{sub:Recontruction}
\item Compute the new momentary proper centroid $\delta\mathbf{d}^{\left(j\right)}$:
\nameref{sub:Proper-centroid-1}; note that the value for $\delta\mathbf{d}$
according to formula (\ref{eq:proper-centroid-1}) is to be assigned
to $\delta\mathbf{d}^{\left(j\right)}$ 
\item Update the position vectors $\mathbf{r}_{k}\longrightarrow\mathbf{r}_{k}-\delta\mathbf{d}^{\left(j\right)}$,
$k=1\ldots n$ of the contour set
\item Calculate the double energy (\ref{eq:minproblem-fine}), exit if the
change (wrt its previous value) is small; Otherwise set $\delta\mathbf{d}^{\left(j+1\right)}=0$,
set $j=j+1$ and repeat from 2.\end{enumerate}
\end{algorithm}

\section{Illustrative examples\label{sec:Illustrative-examples}}

\begin{figure} [!h]
 \centering
 \includegraphics[width=8cm]{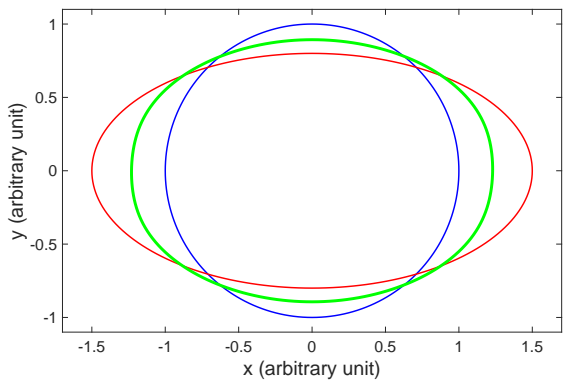}
 \caption{ Mean contour (green) calculated from the RPSV representations of a circle and an ellipse. The mean contour does not pass the intersection of the constituents. }
 \label{fig:representation_set}
\end{figure}

\begin{figure}[!h]
 \centering

\begin{tabular}{c c}
	\includegraphics[width=6cm]{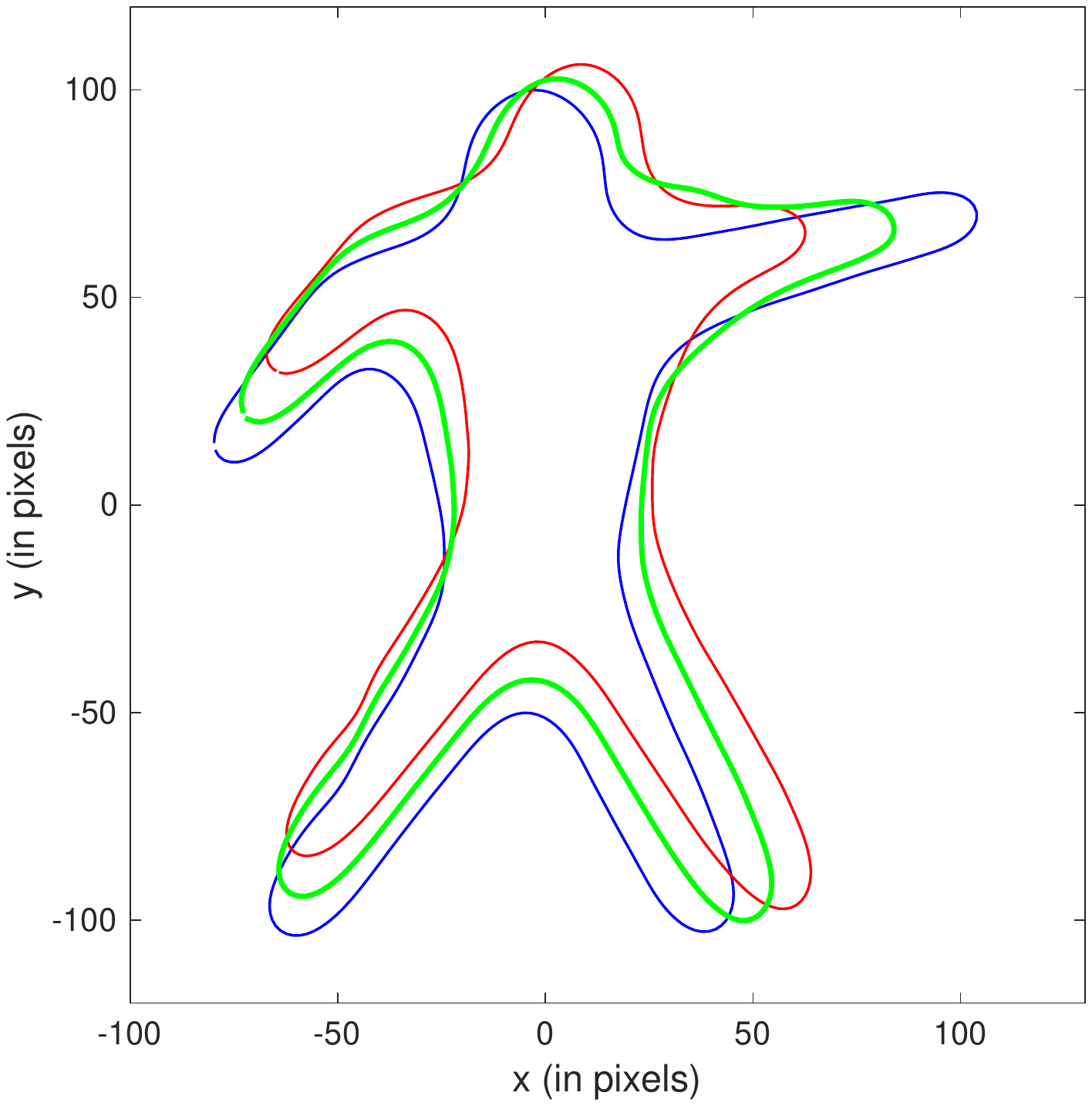}&
	\includegraphics[width=6cm]{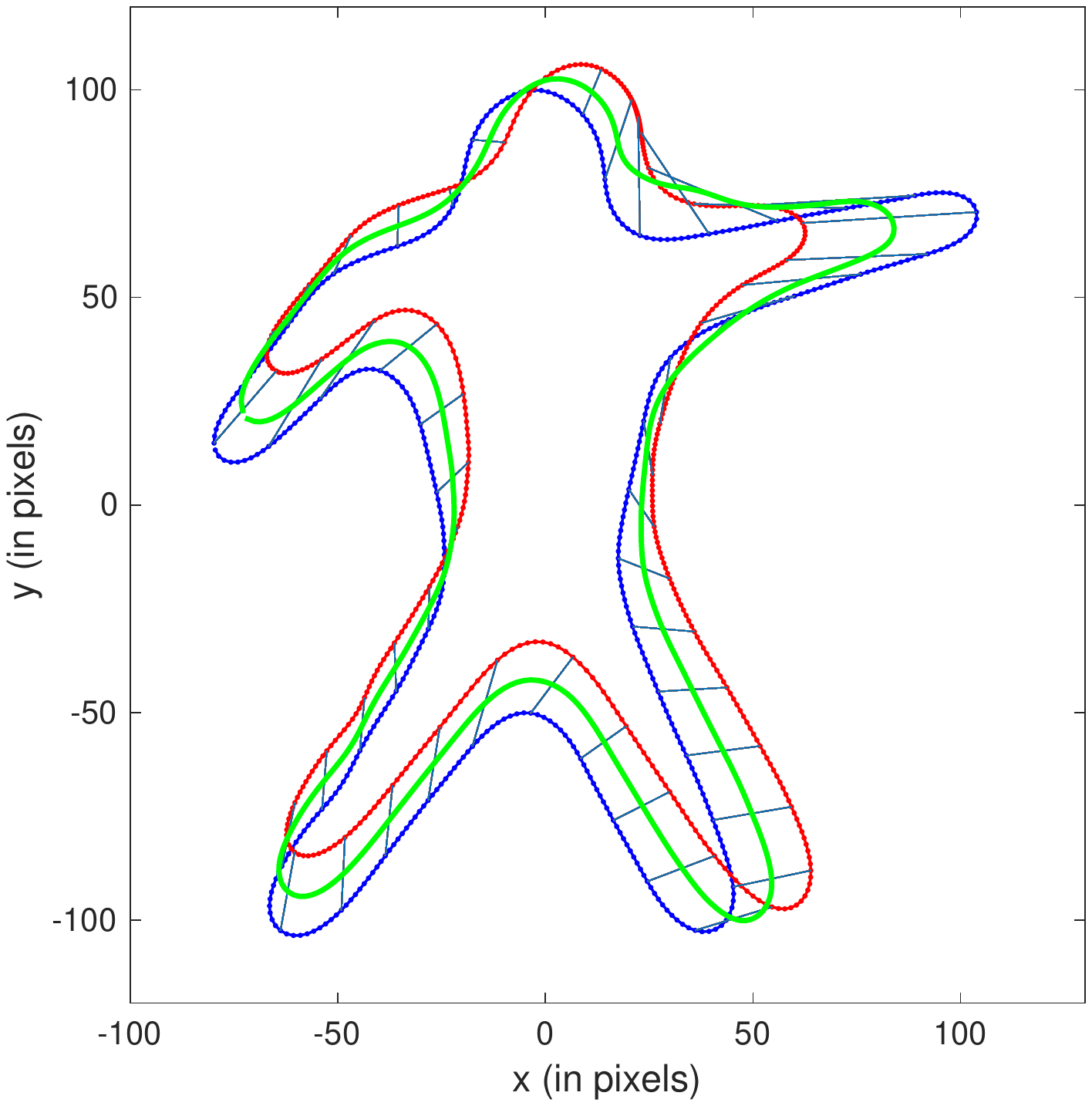}\\
 \end{tabular}

 \caption{Running man without and with marking of point corespondences.}
 \label{fig:test_figure1}
\end{figure}

The illustrations show mean of representation $\left|\mathbf{r}\right|^{m}\mathbf{u}\sqrt{\left|\dot{\mathbf{r}}\right|}$
for $m=1$ \cref{fig:representation_set} a) for one of the simplest
circle/ellipse case (notice that the mean contour does not pass the
intersection of the constituents), b) the mean of non-trivial contours
without and with marking point corespondences \cref{fig:test_figure1}

\section{Conclusion\label{sec:Conclusion}}

In this paper a contour mean determination method - that designed
for averaging manual delineation of objects having non definit boundaries
- was presented. The mean contour is calculated from a set of contours
in a way that all visible information (relative placement, rotation,
scale) are retained. At the same time - borrowed the idea from the
state of the art shape analysis methods - the contour parameterization
is relaxed. The chosen contour representation (RPSV) and the imposed
$\mathbb{L}^{2}$ metric forms a Hilbert space of the contour representations.
The metric is chosen to be invariant wrt the reparameterization, the
distance function based on it has well defined meaning, the (sum of)
the second moment of the contours. The mean contour calculation is
performed in the quotient space space of contours modulo reparameterization
group and could be formulated as a double optimization problem: a
variational for the system of the optimal parameterization and an
extreme value problem for the proper centroid identification. Illustrative
examples show that the resulted mean contours are intuitive according
to human perception sense. Similarities/dissimilarities can be simple
measured and the outlayers determined in this manner are also coincident
with the human perception.

The approach can be generalized in many ways \textit{e.g.} defining
various combination of representations and the associated metrics
(some of them are partly addressed in the article) that may lead meaningful
shape analysis techniques alternative to the current mainstream. Another
plausible direction is the generalization of the method to surfaces.

\part*{Appendices}

In the appendices, the important properties of the action of the reparameterization
group $t\rightarrow\gamma\left(t\right)$, $\mathbf{q}\rightarrow\mathbf{q}\circ\gamma$
\nameref{sec:Appendix-A} and the founding theorems of the mean contour
calculation \nameref{sec:Appendix-B} are examined. The reconstruction
equations are derived in \nameref{sec:Appendix-C}.

Notations and terminology used throughout the appendices are as follow.
Curves are given by their position vectors wrt some standard basis
$\mathbf{i},\,\mathbf{j}$ and denoted as $\mathbf{r}\left(t\right)=x\left(t\right)\mathbf{i}+y\left(t\right)\mathbf{j}$
where $x\left(t\right),\, y\left(t\right)$ are the coordinate functions;
contours are closed curves: $\mathbf{r}\left(0\right)=\mathbf{r}\left(T\right)$.
The discrete representation of a contour is given by the set of $M$
points selected at parameter values distributed uniformly, that is:
$\mathbf{r}_{1}=\mathbf{r}\left(t_{1}\right),\ldots\mathbf{r}_{M}=\mathbf{r}\left(t_{m}\right)$,
$t_{i+1}-t_{i}=\varDelta t$, $t_{1}=0$, $t_{M}=T-\varDelta t$.

Vectors are written with bold letters; vector juxtaposition $\mathbf{a}\mathbf{b}$
indicates direct (dyadic) product, scalar (contraction of a dyad)
and cross products are denoted by dot $\mathbf{a}\cdot\mathbf{b}$
and cross $\mathbf{a}\times\mathbf{b}$ respectively. Derivatives
wrt contour parameter $t$ are denoted by dots: $\mathbf{\dot{r}}\equiv\frac{d\mathbf{r}}{dt}$,
$\ddot{\mathbf{r}}\equiv\frac{d\mathbf{^{2}r}}{dt^{2}}\ldots$ (and
dot is reserved to denote the derivatives wrt $t$); the derivatives
at $\gamma\left(t\right)$ are denoted by primes: $\mathbf{r}^{\prime}\equiv\frac{d\mathbf{r}}{d\gamma}$,
$\mathbf{r}^{\prime\prime}\equiv\frac{d\mathbf{^{2}r}}{d\gamma^{2}}\ldots$.
For the line integrals along a contour (along closed curve), symbol
$\oint$ is used. In the case of iterative methods, the identifiers
of the iteration ('iteration index') are denoted by upper indices
in parentheses \textit{e.g.} the value of the quantity $x$ in the
$k$-th iteration is $x^{\left(k\right)}$.

\section*{Appendix A\label{sec:Appendix-A}}

Property A1: the reparameterization group $\mathbf{q}\rightarrow\mathbf{q}\circ\gamma$
($t\rightarrow\gamma\left(t\right)$) acts by isometries wrt the chosen
representation $\mathbf{q}=\mathbf{r}\sqrt{\left|\mathbf{\dot{r}}\right|}$
($\mathbf{\dot{r}}\equiv\frac{d\mathbf{r}}{dt}$) and metric $d^{2}\left(\mathbf{q}_{1},\mathbf{q}_{2}\right)\doteq\oint\left(\mathbf{r}_{1}\left(t\right)\sqrt{\left|\mathbf{\dot{r}}_{1}\left(t\right)\right|}-\mathbf{r}_{2}\left(t\right)\sqrt{\left|\mathbf{\dot{r}}_{2}\left(t\right)\right|}\right)^{2}dt$.

Proof: consider the common reparameterization $t\rightarrow\gamma\left(t\right)$
of the two contours $\mathbf{q}_{1},\mathbf{q}_{2}$ involved, then
the relation between the operators become $\frac{d}{dt}=\dot{\gamma}\frac{d}{d\gamma}$
($\dot{\gamma}\equiv\frac{d\gamma}{dt}$). The change squared distance
\begin{eqnarray}
\left\Vert \mathbf{q}_{1}\circ\gamma-\mathbf{q}_{2}\circ\gamma\right\Vert ^{2} & = & \oint\left(\mathbf{r}_{1}\left(\gamma\right)\sqrt{\left|\frac{d\mathbf{r}_{1}\left(\gamma\right)}{d\gamma}\right|\dot{\gamma}}-\mathbf{r}_{2}\left(\gamma\right)\sqrt{\left|\frac{d\mathbf{r}_{2}\left(\gamma\right)}{d\gamma}\dot{\gamma}\right|}\right)^{2}dt\nonumber \\
 & = & \oint\left(\mathbf{r}_{1}\left(\gamma\right)\sqrt{\left|\frac{d\mathbf{r}_{1}\left(\gamma\right)}{d\gamma}\right|}-\mathbf{r}_{2}\left(\gamma\right)\sqrt{\left|\frac{d\mathbf{r}_{2}\left(\gamma\right)}{d\gamma}\right|}\right)^{2}\dot{\gamma}dt\\
 & = & \oint\left(\mathbf{r}_{1}\left(\gamma\right)\sqrt{\left|\frac{d\mathbf{r}_{1}\left(\gamma\right)}{d\gamma}\right|}-\mathbf{r}_{2}\left(\gamma\right)\sqrt{\left|\frac{d\mathbf{r}_{2}\left(\gamma\right)}{d\gamma}\right|}\right)^{2}d\gamma\,.\nonumber 
\end{eqnarray}
The last line is equivalent to the definition with renamed variable
of integration, \textit{i.e.} the common reparameterization of the
contours does not influence their distance. This property allows simple
strategy to determine the optimal parameterization system of contours,
that is \textbf{Lemma A1}: one can designate any constituent of the
set of $n$ contours as the reference contour to determine the optimally
parameterized system of contours with pairwise calculation of the
optimal (in the sense of minimum distances) reparameterization wrt
the reference contour.

Proof: assume we have the system of $n$ contours $\mathbf{q}_{1},\,\mathbf{q}_{2},\ldots\mathbf{q}_{n}$
parameterized having same parameter range $\left[0,T\right]$ (otherwise
arbitrarily). First we determine $\gamma_{n}^{\star}$ acting between
$\mathbf{q}_{n-1},\,\mathbf{q}_{n}$ such that $d^{2}\left(\mathbf{q}_{n-1},\mathbf{q}_{n}\circ\gamma_{n}^{\star}\right)$
admits its minimum, second we repeat with $\gamma_{n-1}^{\star}$
such that $d^{2}\left(\mathbf{q}_{n-2},\mathbf{q}_{n-1}\circ\gamma_{n-1}^{\star}\right)$
to be minimal, and update $\mathbf{q}_{n}\circ\gamma_{n}^{\star}\rightarrow\mathbf{q}_{n}\circ\gamma_{n-1}^{\star}\circ\gamma_{n}^{\star}$.
Continuing this procedure, at the end we have the optimally reparameterized
system: $\mathbf{q}_{1},\,\mathbf{q}_{2}\circ\gamma_{2}^{\star},\ldots\mathbf{q}_{n}\circ\gamma_{2}^{\star}\cdots\circ\gamma_{n}^{\star}$.
However, if the pairwise calculations provide unique solution to the
problem $\underset{\gamma_{1i}}{\min}d^{2}\left(\mathbf{q}_{1},\mathbf{q}_{i}\circ\gamma_{i}\right)$,
$i=2,\ldots n$ then the equivalences $\gamma_{1i}\equiv\gamma_{2}^{\star}\cdots\circ\gamma_{i}^{\star}$
must hold. Since the both the reference contour and the order of the
contours are arbitrary, the final system is optimally parameterized
in the minimum distance sense. 

The optimal reparameterization can be uniquely determined, using variational
minimization \textit{e.g.} between contours $1$ and $2$ it can be
formulated as:
\begin{equation}
\underset{\gamma_{12}}{\min}\oint\left(\mathbf{r}_{1}\sqrt{\left|\dot{\mathbf{r}}_{1}\right|}-\mathbf{r}_{2}\left(\gamma_{2}\right)\sqrt{\left|\mathbf{r}_{2}^{\prime}\left(\gamma_{2}\right)\right|\dot{\gamma}_{2}}\right)^{2}dt,\label{eq:optimal-param-statement}
\end{equation}
where the notation $\mathbf{r}_{2}^{\prime}=\frac{d\mathbf{r}_{2}}{d\gamma_{2}}$
is used (dot is exclusively reserved for $t$). The variational problem
is solved via its associated Euler-Lagrange equation.

Property A2: along a linear path $\left(1-\tau\right)\mathbf{q}_{1}+\tau\mathbf{q}_{2}$
the same Euler-Lagrange equation determines the minimal distance solution
between (any) two endpoints $\mathbf{q}_{1},\,\mathbf{q}_{2}$. 

Proof: the distance minimizer integral for the point $\left(1-\tau\right)\mathbf{q}_{1}+\tau\mathbf{q}_{2}$
is: 
\begin{gather}
\oint\left\{ \mathbf{r}_{1}\sqrt{\left|\dot{\mathbf{r}}_{1}\right|}-\left[\left(1-\tau\right)\mathbf{r}_{1}\sqrt{\left|\dot{\mathbf{r}}_{1}\right|}+\tau\mathbf{r}_{2}\left(\gamma_{2}\right)\sqrt{\left|\mathbf{r}_{2}^{\prime}\left(\gamma_{2}\right)\right|\dot{\gamma}_{2}}\right]\right\} ^{2}dt\nonumber \\
\qquad\qquad\qquad\qquad\qquad=\tau^{2}\oint\left(\mathbf{r}_{1}\sqrt{\left|\dot{\mathbf{r}}_{1}\right|}-\mathbf{r}_{2}\left(\gamma_{2}\right)\sqrt{\left|\mathbf{r}_{2}^{\prime}\left(\gamma_{2}\right)\right|\dot{\gamma}_{2}}\right)^{2}dt\label{eq:linear-path-distance}
\end{gather}
the right side differ from the functional to be minimized (\ref{eq:optimal-param-statement})
only in a constant factor which does not affect the associated Euler-Lagrange
equation.

Property A3: also, it is obvious from (\ref{eq:optimal-param-statement})
that the distance ($\sqrt{d^{2}}$) along a linear path alters linearly.

\section*{Appendix B\label{sec:Appendix-B}}

Let $\mathbf{q}_{1},\,\mathbf{q}_{2},\cdots\mathbf{q}_{n}$, $\mathbf{q}_{k}=\mathbf{r}_{k}\sqrt{\dot{\mathbf{r}}_{k}}$
a system of representations of $n$ contours. We wish to determine
the system of optimal reparameterization $\gamma_{k}$, $k=1..n$
that minimizes the squared distances $d^{2}\left(\mathbf{q}_{i},\mathbf{q}_{k}\right)$,
$i,\, k=1,\ldots n$ (\ref{eq:squared-distance}) between them. It
can be done pairwise wrt a reference contour (see \nameref{sec:Appendix-A}).
Without loss of generality, let $\mathbf{r}_{1}$ (represented by
$\mathbf{q}_{1}$) be chosen as the reference contour (hence $\gamma_{1}\left(t\right)\equiv t$),
then the functionals $\oint\left(\mathbf{r}_{1}\sqrt{\left|\dot{\mathbf{r}}_{1}\right|}-\mathbf{r}_{k}\sqrt{\left|\dot{\mathbf{r}}_{k}\right|}\right)^{2}dt$,
$k=2,\ldots n$ are to be minimized wrt the $k$-th diffeomorphism
$\gamma_{k}=\gamma_{k}\left(t\right)$. 

Using the notations (and dependencies on the different contour parameters)
listed below
\begin{eqnarray}
\mathbf{r}_{1} & = & \mathbf{r}_{1}\left(t\right)\nonumber \\
\mathbf{r}_{k} & = & \mathbf{r}_{k}\left(\gamma_{k}\right)\,,\gamma_{k}=\gamma_{k}\left(t\right)\nonumber \\
\mathbf{\dot{r}}_{1} & = & \mathbf{\dot{r}}_{1}\left(t\right)=\frac{d\mathbf{r}_{1}\left(t\right)}{dt}\\
\mathbf{\dot{r}}_{k} & = & \mathbf{\dot{r}}_{k}\left(t\right)=\dot{\gamma}_{k}\left(t\right)\frac{d\mathbf{r}_{k}\left(\gamma_{k}\right)}{d\gamma_{k}}=\dot{\gamma}_{k}\mathbf{r}_{k}^{\prime}\nonumber \\
\mathbf{e}_{k} & = & \frac{\mathbf{\dot{r}}_{k}}{\left|\mathbf{\dot{r}}_{k}\right|}=\frac{\mathbf{r}_{k}^{\prime}}{\left|\mathbf{r}_{k}^{\prime}\right|}\nonumber 
\end{eqnarray}
we first state \textbf{Lemma B1}: The Euler-Lagrange equation associated
with the minimization problem $\underset{\gamma_{k}}{\min}\oint\left(\mathbf{r}_{1}\sqrt{\left|\dot{\mathbf{r}}_{1}\right|}-\mathbf{r}_{k}\left(\gamma_{k}\right)\sqrt{\left|\mathbf{r}_{k}^{\prime}\left(\gamma_{k}\right)\right|\dot{\gamma}_{k}}\right)^{2}dt$
is $\dot{\mathbf{r}}_{1}\cdot\mathbf{r}_{k}-\mathbf{r}_{1}\cdot\dot{\mathbf{r}}_{k}+\frac{1}{2}\mathbf{r}_{1}\cdot\mathbf{r}_{k}\left(\varGamma_{1}-\varGamma_{k}\right)$.

Proof: the Lagranian and its derivatives are:
\begin{eqnarray}
L\left(\gamma_{k},\dot{\gamma}_{k}\right) & = & \left(\mathbf{r}_{1}\sqrt{\left|\dot{\mathbf{r}}_{1}\right|}-\mathbf{r}_{k}\sqrt{\left|\mathbf{r}_{k}^{\prime}\right|\dot{\gamma}_{k}}\right)^{2}\nonumber \\
\frac{\partial L}{\partial\gamma_{k}} & =- & 2\left(\mathbf{r}_{1}\sqrt{\left|\dot{\mathbf{r}}_{1}\right|}-\mathbf{r}_{k}\sqrt{\left|\mathbf{r}_{k}^{\prime}\right|\dot{\gamma}_{k}}\right)\cdot\left(\mathbf{r}_{k}^{\prime}\sqrt{\left|\mathbf{r}_{k}^{\prime}\right|\dot{\gamma}_{k}}+\mathbf{r}_{k}\frac{\dot{\gamma}_{k}\mathbf{e}_{k}\cdot\mathbf{r}_{k}^{\prime\prime}}{2\sqrt{\left|\mathbf{r}_{k}^{\prime}\right|\dot{\gamma}_{k}}}\right)\nonumber \\
\frac{\partial L}{\partial\dot{\gamma}_{k}} & = & -\left(\mathbf{r}_{1}\sqrt{\left|\dot{\mathbf{r}}_{1}\right|}-\mathbf{r}_{k}\sqrt{\left|\mathbf{r}_{k}^{\prime}\right|\dot{\gamma}_{k}}\right)\cdot\mathbf{r}_{k}\frac{\left|\mathbf{r}_{k}^{\prime}\right|}{\sqrt{\left|\mathbf{r}_{k}^{\prime}\right|\dot{\gamma}_{k}}}
\end{eqnarray}
From the relations between the differential operators
\begin{eqnarray}
\frac{d}{d\gamma} & = & \frac{1}{\dot{\gamma}}\frac{d}{dt}\nonumber \\
\frac{d^{2}}{d\gamma^{2}} & = & \frac{1}{\dot{\gamma}}\left(-\frac{\ddot{\gamma}}{\dot{\gamma}^{2}}\frac{d}{dt}+\frac{1}{\dot{\gamma}}\frac{d^{2}}{dt^{2}}\right),\label{eq:operator-relation}
\end{eqnarray}
we have
\begin{eqnarray*}
\frac{\partial L}{\partial\gamma_{k}} & = & -2\left(\mathbf{r}_{1}\sqrt{\left|\dot{\mathbf{r}}_{1}\right|}-\mathbf{r}_{k}\sqrt{\left|\dot{\mathbf{r}}_{k}\right|}\right)\cdot\left[\frac{1}{\dot{\gamma}_{k}}\dot{\mathbf{r}}_{k}\sqrt{\left|\dot{\mathbf{r}}_{k}\right|}+\mathbf{r}_{k}\frac{\mathbf{e}_{k}\cdot\left(-\frac{\ddot{\gamma}}{\dot{\gamma}^{2}}\dot{\mathbf{r}}_{k}+\frac{1}{\dot{\gamma}}\ddot{\mathbf{r}}_{k}\right)}{2\sqrt{\left|\dot{\mathbf{r}}_{k}\right|}}\right]\\
 & = & -\left(\mathbf{r}_{1}\sqrt{\left|\dot{\mathbf{r}}_{1}\right|}-\mathbf{r}_{k}\sqrt{\left|\dot{\mathbf{r}}_{k}\right|}\right)\cdot\frac{\sqrt{\left|\dot{\mathbf{r}}_{k}\right|}}{\dot{\gamma}_{k}}\left(2\dot{\mathbf{r}}_{k}-\frac{\ddot{\gamma}}{\dot{\gamma}}\mathbf{r}_{k}+\mathbf{r}_{k}\frac{\dot{\mathbf{r}}_{k}\cdot\ddot{\mathbf{r}}_{k}}{\left|\dot{\mathbf{r}}_{k}\right|^{2}}\right)\\
 & = & -\frac{1}{\dot{\gamma}_{k}}\left(\mathbf{r}_{1}\sqrt{\left|\dot{\mathbf{r}}_{1}\right|\left|\dot{\mathbf{r}}_{k}\right|}-\mathbf{r}_{k}\left|\dot{\mathbf{r}}_{k}\right|\right)\cdot\left(2\dot{\mathbf{r}}_{k}-\frac{\ddot{\gamma}}{\dot{\gamma}}\mathbf{r}_{k}+\mathbf{r}_{k}\frac{\dot{\mathbf{r}}_{k}\cdot\ddot{\mathbf{r}}_{k}}{\left|\dot{\mathbf{r}}_{k}\right|^{2}}\right)\\
\frac{\partial L}{\partial\dot{\gamma}_{k}} & = & -\left(\mathbf{r}_{1}\sqrt{\left|\dot{\mathbf{r}}_{1}\right|}-\mathbf{r}_{k}\sqrt{\left|\dot{\mathbf{r}}_{k}\right|}\right)\cdot\mathbf{r}_{k}\frac{\sqrt{\left|\dot{\mathbf{r}}_{k}\right|}}{\dot{\gamma}_{k}}\\
 & = & -\frac{1}{\dot{\gamma}_{k}}\left(\mathbf{r}_{1}\sqrt{\left|\dot{\mathbf{r}}_{1}\right|\left|\dot{\mathbf{r}}_{k}\right|}-\mathbf{r}_{k}\left|\dot{\mathbf{r}}_{k}\right|\right)\cdot\mathbf{r}_{k},
\end{eqnarray*}
and
\begin{eqnarray*}
-\frac{d}{dt}\frac{\partial L}{\partial\dot{\gamma}_{k}} & = & \left(\mathbf{r}_{1}\sqrt{\left|\dot{\mathbf{r}}_{1}\right|\left|\dot{\mathbf{r}}_{k}\right|}-\mathbf{r}_{k}\left|\dot{\mathbf{r}}_{k}\right|\right)\cdot\left(\frac{1}{\dot{\gamma}_{k}}\mathbf{\dot{r}}_{k}-\frac{\ddot{\gamma}}{\dot{\gamma}_{k}^{2}}\mathbf{r}_{k}\right)\\
 &  & \qquad+\frac{\sqrt{\left|\dot{\mathbf{r}}_{1}\right|\left|\dot{\mathbf{r}}_{k}\right|}}{\dot{\gamma}_{k}}\mathbf{r}_{k}\cdot\left(\dot{\mathbf{r}}_{1}+\frac{1}{2}\mathbf{r}_{1}\frac{\mathbf{\dot{r}}_{1}\cdot\ddot{\mathbf{r}}_{1}}{\left|\dot{\mathbf{r}}_{1}\right|^{2}}+\frac{1}{2}\mathbf{r}_{1}\frac{\mathbf{\dot{r}}_{k}\cdot\ddot{\mathbf{r}}_{k}}{\left|\dot{\mathbf{r}}_{k}\right|^{2}}\right)\\
 &  & \qquad\qquad-\frac{1}{\dot{\gamma}_{k}}\mathbf{r}_{k}\cdot\left(\mathbf{\dot{r}}_{k}\left|\dot{\mathbf{r}}_{k}\right|+\mathbf{r}_{k}\frac{\dot{\mathbf{r}}_{k}\cdot\ddot{\mathbf{r}}_{k}}{\left|\dot{\mathbf{r}}_{k}\right|^{2}}\right)\,.
\end{eqnarray*}
The Euler-Lagrange equation for the k-th diffeomorphism $\gamma_{k}=\gamma_{k}\left(t\right)$
is:
\begin{eqnarray*}
\frac{\partial L}{\partial\dot{\gamma}_{k}}-\frac{d}{dt}\frac{\partial L}{\partial\dot{\gamma}_{k}} & = & \frac{\sqrt{\left|\dot{\mathbf{r}}_{1}\right|\left|\dot{\mathbf{r}}_{k}\right|}}{\dot{\gamma}_{k}}\left[\dot{\mathbf{r}}_{1}\cdot\mathbf{r}_{k}-\mathbf{r}_{1}\cdot\dot{\mathbf{r}}_{k}+\frac{1}{2}\mathbf{r}_{1}\cdot\mathbf{r}_{k}\left(\frac{\mathbf{\dot{r}}_{1}\cdot\ddot{\mathbf{r}}_{1}}{\left|\dot{\mathbf{r}}_{1}\right|^{2}}-\frac{\mathbf{\dot{r}}_{k}\cdot\ddot{\mathbf{r}}_{k}}{\left|\dot{\mathbf{r}}_{k}\right|^{2}}\right)\right]\,.
\end{eqnarray*}
Assuming $\mathbf{\frac{\sqrt{\left|\dot{\mathbf{r}}_{1}\right|\left|\dot{\mathbf{r}}_{k}\right|}}{\dot{\gamma}_{k}}}$
is not zero at any point, we can divide with it, then the EulerLagrange
equations to be solved are given with:
\begin{eqnarray}
\dot{\mathbf{r}}_{1}\cdot\mathbf{r}_{k}-\mathbf{r}_{1}\cdot\dot{\mathbf{r}}_{k}+\frac{1}{2}\mathbf{r}_{1}\cdot\mathbf{r}_{k}\left(\varGamma_{1}-\varGamma_{k}\right) & = & 0,\, k=2,\ldots n,\label{eq:Euler-Lagrange}
\end{eqnarray}
where 'Christoffel divergences' $\varGamma_{i}=\frac{\mathbf{\dot{r}}_{i}\cdot\ddot{\mathbf{r}}_{i}}{\left|\dot{\mathbf{r}}_{i}\right|^{2}}$,
$i=1,\ldots n$ are introduced to simplify the equation. 

Note that the optimal contour system can be generalized in many ways,
\textit{e.g.} for the representation $\mathbf{q}=f\left(\mathbf{r}\right)\mathbf{u}\sqrt{\left|\dot{\mathbf{r}}\right|}$
- where $\mathbf{u}=\frac{\mathbf{r}}{\left|\mathbf{r}\right|}$ is
the unit vector in the direction of the position vector, $f$ is appropriately
defined scalar valued function. Here we provide equations for the
$\mathbf{q}=\left|\mathbf{r}\right|^{m}\mathbf{u}\sqrt{\left|\dot{\mathbf{r}}\right|}$,
$m\in\mathbb{R}$ cases (the $m=1\,\rightarrow$ $\mathbf{q}=\mathbf{r}\sqrt{\left|\dot{\mathbf{r}}\right|}$
is the case examined in this paper in details). For these cases, the
pairwise distance minimizers based on the $\mathbb{L}^{2}$ metric
are formulated as: 
\[
\underset{\gamma_{k}}{\min}\oint\left(\left|\mathbf{r}_{1}\right|^{m}\mathbf{u}_{1}\sqrt{\left|\dot{\mathbf{r}}_{1}\right|}-\left|\mathbf{r}_{k}\left(\gamma_{k}\right)\right|^{m}\mathbf{u}_{k}\left(\gamma_{k}\right)\sqrt{\left|\mathbf{r}_{k}^{\prime}\left(\gamma_{k}\right)\right|\dot{\gamma}_{k}}\right)^{2}dt,
\]
and the associated Euler-Lagrange equations take the form:
\begin{equation}
\dot{\mathbf{r}}_{1}\cdot\left(m\mathbf{u}_{1}\mathbf{u}_{1}+\mathbf{u}_{1}^{\perp}\mathbf{u}_{1}^{\perp}\right)\cdot\mathbf{r}_{k}-\mathbf{r}_{1}\cdot\left(m\mathbf{u}_{k}\mathbf{u}_{k}+\mathbf{u}_{k}^{\perp}\mathbf{u}_{k}^{\perp}\right)\cdot\dot{\mathbf{r}}_{k}+\frac{1}{2}\mathbf{r}_{1}\cdot\mathbf{r}_{k}\left(\varGamma_{1}-\varGamma_{k}\right)=0
\end{equation}
where $\mathbf{u}^{\perp}=\mathbf{k}\times\mathbf{u}$ is the unit
vector perpendicular to the position vector ($\mathbf{k}$ is the
unit normal of the plane). There is singularity at $m=-\frac{1}{2}$
(a uniform scaling $\mathbf{r}\rightarrow\alpha\mathbf{r}$ leads
to the same representation $\left.\mathbf{q}\right|_{\alpha\mathbf{r}}=\left.\mathbf{q}\right|_{\mathbf{r}}=\mathbf{u}\sqrt{\frac{\left|\dot{\mathbf{r}}\right|}{\left|\mathbf{r}\right|}}$).
For this value the reconstruction cannot be made (see also \nameref{sec:Appendix-C}).

The important consequence of the Lemma B1: 

\textbf{Theorem B2}: the solution for the minimization problem (\ref{eq:minproblem-raw})
$\underset{\gamma_{i}}{\min}\stackrel[i=1]{n}{\sum}\oint\left(\mathbf{r}\left(t\right)\sqrt{\left|\dot{\mathbf{r}}\left(t\right)\right|}-\mathbf{r}_{i}\left(t_{i}\right)\sqrt{\left|\dot{\mathbf{r}}_{i}\left(t_{i}\right)\right|}\right)^{2}dt$,
where $\mathbf{r}\left(t\right)\sqrt{\left|\dot{\mathbf{r}}\left(t\right)\right|}=\frac{1}{n}\stackrel[i=1]{n}{\sum}\mathbf{r}_{i}\left(t_{i}\right)\sqrt{\left|\dot{\mathbf{r}}_{i}\left(t_{i}\right)\right|}$
is the system of optimal reparameterization $t_{i}=\gamma_{i}\left(t\right)$,
$i=1,\ldots n$ determined by the pairwise optimizations between the
constituents.

Proof: a) repeating the steps of the previous proof, the optimal parameterization
system satisfies the set of Euler-Lagrange equations: 

\begin{eqnarray}
\dot{\mathbf{r}}\cdot\mathbf{r}_{k}-\mathbf{r}\cdot\dot{\mathbf{r}}_{k}+\frac{1}{2}\mathbf{r}\cdot\mathbf{r}_{k}\left(\frac{\mathbf{\dot{r}}\cdot\ddot{\mathbf{r}}}{\left|\dot{\mathbf{r}}\right|^{2}}-\frac{\mathbf{\dot{r}}_{k}\cdot\ddot{\mathbf{r}}_{k}}{\left|\dot{\mathbf{r}}_{k}\right|^{2}}\right) & = & 0,\, k=1,\ldots n,\label{eq:Euler-Lagrange-1}
\end{eqnarray}
b) taking the derivative wrt $t$ of the mean expression $\mathbf{r}\left(t\right)\sqrt{\left|\dot{\mathbf{r}}\left(t\right)\right|}=\frac{1}{n}\stackrel[i=1]{n}{\sum}\mathbf{r}_{i}\left(t_{i}\right)\sqrt{\left|\dot{\mathbf{r}}_{i}\left(t_{i}\right)\right|}$
then the dot product with $\mathbf{r}_{k}$, we have:
\begin{eqnarray}
\sqrt{\mathbf{\left|\dot{\mathbf{r}}\right|}}\left(\dot{\mathbf{r}}\cdot\mathbf{r}_{k}+\frac{1}{2}\mathbf{r}\cdot\mathbf{r}_{k}\frac{\mathbf{\dot{r}}\cdot\ddot{\mathbf{r}}}{\left|\dot{\mathbf{r}}\right|^{2}}\right) & = & \frac{1}{n}\stackrel[i=1]{n}{\sum}\sqrt{\left|\mathbf{\dot{\mathbf{r}}}_{i}\right|}\left(\dot{\mathbf{r}}_{i}\cdot\mathbf{r}_{k}+\frac{1}{2}\mathbf{r}_{i}\cdot\mathbf{r}_{k}\frac{\mathbf{\dot{r}}_{i}\cdot\ddot{\mathbf{r}}_{i}}{\left|\dot{\mathbf{r}}_{i}\right|^{2}}\right)\,.\label{eq:partial-result-1}
\end{eqnarray}
As assumed (\ref{eq:Euler-Lagrange}) equations are satisfied. From
this 

\begin{eqnarray}
\dot{\mathbf{r}}_{i}\cdot\mathbf{r}_{k}+\frac{1}{2}\mathbf{r}_{i}\cdot\mathbf{r}_{k}\frac{\mathbf{\dot{r}}_{i}\cdot\ddot{\mathbf{r}}_{i}}{\left|\dot{\mathbf{r}}_{i}\right|^{2}} & = & \mathbf{r}_{i}\cdot\dot{\mathbf{r}}_{k}+\frac{1}{2}\mathbf{r}_{i}\cdot\mathbf{r}_{k}\frac{\mathbf{\dot{r}}_{k}\cdot\ddot{\mathbf{r}}_{k}}{\left|\dot{\mathbf{r}}_{k}\right|^{2}}\,.\label{eq:partial-result-2}
\end{eqnarray}
Substituting (\ref{eq:partial-result-2}) to (\ref{eq:partial-result-1}),
we get:
\begin{gather}
\sqrt{\mathbf{\left|\dot{\mathbf{r}}\right|}}\left(\dot{\mathbf{r}}\cdot\mathbf{r}_{k}+\frac{1}{2}\mathbf{r}\cdot\mathbf{r}_{k}\frac{\mathbf{\dot{r}}\cdot\ddot{\mathbf{r}}}{\left|\dot{\mathbf{r}}\right|^{2}}\right)=\frac{1}{n}\stackrel[i=1]{n}{\sum}\sqrt{\left|\mathbf{\dot{\mathbf{r}}}_{i}\right|}\left(\mathbf{r}_{i}\cdot\dot{\mathbf{r}}_{k}+\frac{1}{2}\mathbf{r}_{i}\cdot\mathbf{r}_{k}\frac{\mathbf{\dot{r}}_{k}\cdot\ddot{\mathbf{r}}_{k}}{\left|\dot{\mathbf{r}}_{k}\right|^{2}}\right)\nonumber \\
\qquad\qquad\qquad=\dot{\mathbf{r}}_{k}\cdot\left(\frac{1}{n}\stackrel[i=1]{n}{\sum}\sqrt{\left|\mathbf{\dot{\mathbf{r}}}_{i}\right|}\mathbf{r}_{i}\right)+\frac{1}{2}\frac{\mathbf{\dot{r}}_{k}\cdot\ddot{\mathbf{r}}_{k}}{\left|\dot{\mathbf{r}}_{k}\right|^{2}}\mathbf{r}_{k}\cdot\left(\frac{1}{n}\stackrel[i=1]{n}{\sum}\sqrt{\left|\mathbf{\dot{\mathbf{r}}}_{i}\right|}\mathbf{r}_{i}\right)\nonumber \\
\qquad\qquad\qquad\qquad=\dot{\mathbf{r}}_{k}\cdot\mathbf{r}\sqrt{\left|\dot{\mathbf{r}}\right|}+\frac{1}{2}\frac{\mathbf{\dot{r}}_{k}\cdot\ddot{\mathbf{r}}_{k}}{\left|\dot{\mathbf{r}}_{k}\right|^{2}}\mathbf{r}_{k}\cdot\mathbf{r}\sqrt{\left|\dot{\mathbf{r}}\right|}\,.
\end{gather}
Rearranging, we have:
\begin{eqnarray*}
\sqrt{\mathbf{\left|\dot{\mathbf{r}}\right|}}\left[\dot{\mathbf{r}}\cdot\mathbf{r}_{k}-\dot{\mathbf{r}}_{k}\cdot\mathbf{r}+\frac{1}{2}\mathbf{r}\cdot\mathbf{r}_{k}\left(\frac{\mathbf{\dot{r}}\cdot\ddot{\mathbf{r}}}{\left|\dot{\mathbf{r}}\right|^{2}}-\frac{\mathbf{\dot{r}}_{k}\cdot\ddot{\mathbf{r}}_{k}}{\left|\dot{\mathbf{r}}_{k}\right|^{2}}\right)\right] & = & 0,
\end{eqnarray*}
that is the $k$-th equation of (\ref{eq:Euler-Lagrange-1}).

\section*{Appendix C\label{sec:Appendix-C}}

In this section we derive the equations used to reconstruct the contours
from their RPSV representation $\mathbf{q}\left(t\right)\rightarrow\mathbf{r}\left(t\right)$,
where $\mathbf{q}\left(t\right)=\mathbf{r}\left(t\right)\sqrt{\left|\mathbf{\dot{r}}\left(t\right)\right|}$
is known. Observing that $\frac{\mathbf{q}\left(t\right)}{\left|\mathbf{q}\left(t\right)\right|}=\frac{\mathbf{r}\left(t\right)}{\left|\mathbf{r}\left(t\right)\right|}$,
we introduce the notation for the unit vector pointing from the proper
centroid to the direction of both points $\mathbf{q}\left(t\right)$,
$\mathbf{r}\left(t\right)$: 
\begin{equation}
\mathbf{u}\left(t\right)\doteq\frac{\mathbf{q}\left(t\right)}{\left|\mathbf{q}\left(t\right)\right|}=\frac{\mathbf{r}\left(t\right)}{\left|\mathbf{r}\left(t\right)\right|}\,.
\end{equation}
Having the direction of the position vector, we need to determine
only its distance measured from the centroid $\left|\mathbf{r}\left(t\right)\right|$
then position vector $\mathbf{r}\left(t\right)=\left|\mathbf{r}\left(t\right)\right|\mathbf{u}\left(t\right)$.%
\footnote{This also means that the unit direction vector remains always constant
(i.e. does not change during the iteration described in this appendix).%
} (Hereinafter we will also use the notation $\mathbf{e}\left(t\right)\doteq\frac{\mathbf{\dot{r}}\left(t\right)}{\left|\mathbf{\dot{r}}\left(t\right)\right|}$
for the unit tangent vector of the contour.) Now we define the scalar
function
\begin{equation}
f\left(\mathbf{r},\dot{\mathbf{r}}\right)\doteq\left|\mathbf{q}\right|-\left|\mathbf{r}\right|\sqrt{\left|\mathbf{\dot{r}}\right|}\,.\label{eq:cost-for-reconstruct}
\end{equation}
With this definition, the determination of $\left|\mathbf{r}\right|$
becomes root finding problem (at each parameter value $t$). In function
(\ref{eq:cost-for-reconstruct}) temporarily we handle the position
vector $\mathbf{r}$ and its derivative $\dot{\mathbf{r}}$ as if
they were independent variables.

Assume we know the value of $f$ at some initial guess point $\mathbf{r}^{\left(k\right)}$,
$\dot{\mathbf{r}}^{\left(k\right)}$ close to its root, its linear
approximation around can be written as
\begin{equation}
y\left(\mathbf{r},\dot{\mathbf{r}}\right)=f+\frac{\partial f}{\partial\mathbf{r}}\cdot\left(\mathbf{r}-\mathbf{r}^{\left(k\right)}\right)+\frac{\partial f}{\partial\dot{\mathbf{r}}}\cdot\left(\mathbf{\dot{r}}-\mathbf{\dot{r}}^{\left(k\right)}\right)\label{eq:cost-linear-approximation}
\end{equation}
where function $f$ and its gradients $\frac{\partial f}{\partial\mathbf{r}}$
and $\frac{\partial f}{\partial\dot{\mathbf{r}}}$ are all evaluated
at $\mathbf{r}^{\left(k\right)}$, $\dot{\mathbf{r}}^{\left(k\right)}$.
The gradients are:
\begin{eqnarray}
\frac{\partial f}{\partial\mathbf{r}} & = & -\sqrt{\left|\mathbf{\dot{r}}\right|}\mathbf{u}\nonumber \\
\frac{\partial f}{\partial\dot{\mathbf{r}}} & = & -\frac{1}{2}\frac{\left|\mathbf{r}\right|}{\sqrt{\left|\mathbf{\dot{r}}\right|}}\mathbf{e}\,.\label{eq:cost-gradients}
\end{eqnarray}
Substituting the gradient expressions into (\ref{eq:cost-linear-approximation})
at point $\mathbf{r}^{\left(k\right)}$, $\dot{\mathbf{r}}^{\left(k\right)}$,
we have the the equation for the root ($y=0$) of the linear approximation
(\ref{eq:cost-linear-approximation}):
\begin{equation}
\left|\mathbf{q}\right|-\left|\mathbf{r}^{\left(k\right)}\right|\sqrt{\left|\mathbf{\dot{r}}^{\left(k\right)}\right|}-\sqrt{\left|\mathbf{\dot{r}}^{\left(k\right)}\right|}\mathbf{u}\cdot\left(\mathbf{r}-\mathbf{r}^{\left(k\right)}\right)-\frac{1}{2}\frac{\left|\mathbf{r}^{\left(k\right)}\right|}{\sqrt{\left|\mathbf{\dot{r}}^{\left(k\right)}\right|}}\mathbf{e}^{\left(k\right)}\cdot\left(\mathbf{\dot{r}}-\mathbf{\dot{r}}^{\left(k\right)}\right)=0\label{eq:linear-cost-to-solve}
\end{equation}
to be solved for $\left|\mathbf{r}\right|$. Using the identities
$\mathbf{e}^{\left(k\right)}\cdot\mathbf{\dot{r}}^{\left(k\right)}\equiv\left|\mathbf{\dot{r}}^{\left(k\right)}\right|$,
$\mathbf{u}\cdot\mathbf{r}^{\left(k\right)}\equiv\left|\mathbf{r}^{\left(k\right)}\right|$,
$\mathbf{u}\cdot\mathbf{r}\equiv\left|\mathbf{r}\right|$, equation
(\ref{eq:linear-cost-to-solve}) can be rearranged as
\begin{equation}
\left|\mathbf{r}\right|+\frac{1}{2}\frac{\left|\mathbf{r}^{\left(k\right)}\right|}{\left|\mathbf{\dot{r}}^{\left(k\right)}\right|}\mathbf{e}^{\left(k\right)}\cdot\mathbf{\dot{r}}=\frac{\left|\mathbf{q}\right|}{\sqrt{\left|\mathbf{\dot{r}}^{\left(k\right)}\right|}}+\frac{1}{2}\left|\mathbf{r}^{\left(k\right)}\right|\,.\label{eq:linear-cost-solution}
\end{equation}
Now we take into account that $\dot{\mathbf{r}}$ is not independent
of $\mathbf{r}$. Assuming our contour (its approximation) is defined
by a discrete set of $M$ points: $\mathbf{r}_{1}=\mathbf{r}\left(t_{1}\right),\ldots\mathbf{r}_{M}=\mathbf{r}\left(t_{M}\right)$,
$t_{1}=0$, $t_{M}=T-\varDelta t$, uniformly distributed wrt $t$,%
\footnote{This assumption is taken throughout the paper.%
} we can introduce the notations for the immediate neighbours of $\mathbf{u}$
and $\mathbf{r}$ at any parameter value $t$ as 
\begin{eqnarray}
\mathbf{u}_{+} & \doteq & \mathbf{u}\left(t+\varDelta t\right)\nonumber \\
\mathbf{u}_{-} & \doteq & \mathbf{u}\left(t-\varDelta t\right)\nonumber \\
\mathbf{r}_{+} & \doteq & \mathbf{r}\left(t+\varDelta t\right)\\
\mathbf{r}_{-} & \doteq & \mathbf{r}\left(t-\varDelta t\right)\,.\nonumber 
\end{eqnarray}
Solution (\ref{eq:linear-cost-solution}) can be approximated using
the simple finite central differences scheme $\mathbf{\dot{r}}\approx\frac{\mathbf{r}_{+}-\mathbf{r}_{-}}{2\varDelta t}$,
$\mathbf{r}_{\pm}=\left|\mathbf{r}_{\pm}\right|\mathbf{u}_{\pm}$
as:
\begin{equation}
\left|\mathbf{r}\right|+\left\{ \frac{1}{4\varDelta t}\frac{\left|\mathbf{r}^{\left(k\right)}\right|}{\left|\mathbf{\dot{r}}^{\left(k\right)}\right|}\mathbf{e}^{\left(k\right)}\cdot\mathbf{u}_{+}\right\} \left|\mathbf{r}_{+}\right|-\left\{ \frac{1}{4\varDelta t}\frac{\left|\mathbf{r}^{\left(k\right)}\right|}{\left|\mathbf{\dot{r}}^{\left(k\right)}\right|}\mathbf{e}^{\left(k\right)}\cdot\mathbf{u}_{-}\right\} \left|\mathbf{r}_{-}\right|=\frac{\left|\mathbf{q}\right|}{\sqrt{\left|\mathbf{\dot{r}}^{\left(k\right)}\right|}}+\frac{1}{2}\left|\mathbf{r}^{\left(k\right)}\right|\,.\label{eq:linear-cost-solution-1}
\end{equation}
On the right side all quantities are known; on the left side the known
coefficients are emphasized by putting them into braces. For the whole
point set this constitutes a linear equation system with sparse matrix
three-diagonal almost everywhere except the first and last line. The
derivation above follows the steps of the derivation of Newton\textendash Raphson
method. This method is widely used to determine the root of the nonlinear
equations iteratively. Starting from an intermediate result (approximation
of the root of (\ref{eq:cost-for-reconstruct})) $\left|\mathbf{r}^{\left(k\right)}\right|$,
the next (expectably more accurate) approximation $\left|\mathbf{r}^{\left(k+1\right)}\right|$
is given as the solution of (\ref{eq:linear-cost-solution-1}). With
the substitution $\left|\mathbf{r}\right|\rightarrow\left|\mathbf{r}^{\left(k+1\right)}\right|$,
the linear equation system $\mathbf{A}^{\left(k\right)}\mathbf{x}^{\left(k+1\right)}=\mathbf{b}^{\left(k\right)}$
needs to be solved for the next ($k+1$-th) root vector $\mathbf{x}^{\left(k+1\right)}$
with the sought components 
\begin{equation}
\mathbf{x}^{\left(k+1\right)}=\left[\begin{array}{ccccc}
\left|\mathbf{r}_{1}^{\left(k+1\right)}\right| & \cdots & \left|\mathbf{r}_{i}^{\left(k+1\right)}\right| & \cdots & \left|\mathbf{r}_{M}^{\left(k+1\right)}\right|\end{array}\right]^{T}\label{eq:ray-length}
\end{equation}
using the matrix 
\begin{equation}
\mathbf{A}^{\left(k\right)}=\left[\begin{array}{ccccc}
1 & a_{1,2}^{\left(k\right)} & 0\cdots & \cdots0 & a_{1,M}^{\left(k\right)}\\
\cdots & \cdots & \cdots & \cdots & \cdots\\
0\cdots & a_{i,i-1}^{\left(k\right)} & 1 & a_{i,i+1}^{\left(k\right)} & \cdots0\\
\cdots & \cdots & \cdots & \cdots & \cdots\\
a_{M,1}^{\left(k\right)} & 0\cdots & \cdots0 & a_{M,M-1}^{\left(k\right)} & 1
\end{array}\right]\label{eq:coeffitient-matrix}
\end{equation}
$a_{i,i-1}^{\left(k\right)}=-\frac{1}{4\varDelta t}\frac{\left|\mathbf{r}_{i}^{\left(k\right)}\right|}{\left|\mathbf{\dot{r}}_{i}^{\left(k\right)}\right|}\mathbf{e}_{i}^{\left(k\right)}\cdot\mathbf{u}_{i-1}$,
$a_{i,i+1}^{\left(k\right)}=\frac{1}{4\varDelta t}\frac{\left|\mathbf{r}_{i}^{\left(k\right)}\right|}{\left|\mathbf{\dot{r}}_{i}^{\left(k\right)}\right|}\mathbf{e}_{i}^{\left(k\right)}\cdot\mathbf{u}_{i+1}$,
$i=2\ldots M-1$, $a_{1,M}^{\left(k\right)}=-\frac{1}{4\varDelta t}\frac{\left|\mathbf{r}_{1}^{\left(k\right)}\right|}{\left|\mathbf{\dot{r}}_{1}^{\left(k\right)}\right|}\mathbf{e}_{1}^{\left(k\right)}\cdot\mathbf{u}_{M}$,
$a_{M,1}^{\left(k\right)}=\frac{1}{4\varDelta t}\frac{\left|\mathbf{r}_{M}^{\left(k\right)}\right|}{\left|\mathbf{\dot{r}}_{M}^{\left(k\right)}\right|}\mathbf{e}_{M}^{\left(k\right)}\cdot\mathbf{u}_{1}$
and the vector
\begin{equation}
\mathbf{b}^{\left(k+1\right)}=\left[\begin{array}{ccccc}
b_{1}^{\left(k\right)} & \cdots & b_{i}^{\left(k\right)} & \cdots & b_{M}^{\left(k\right)}\end{array}\right]^{T},\label{eq:constant-vector}
\end{equation}
$b_{i}^{\left(k\right)}=\frac{\left|\mathbf{q}_{i}\right|}{\sqrt{\left|\mathbf{\dot{r}}_{i}^{\left(k\right)}\right|}}+\frac{1}{2}\left|\mathbf{r}_{i}^{\left(k\right)}\right|$,
$i=1\ldots M$ calculable from the $k$-th iteration.

Note that for the generalized representation $\mathbf{q}=\left|\mathbf{r}\right|^{m}\mathbf{u}\sqrt{\left|\dot{\mathbf{r}}\right|}$
the reconstruction equations (\ref{eq:linear-cost-solution-1}) (with
the substitution $\left|\mathbf{r}\right|\rightarrow\left|\mathbf{r}^{\left(k+1\right)}\right|$)
take the form:
\begin{gather}
\left|\mathbf{r}^{\left(k+1\right)}\right|+\left\{ \frac{1}{4m\varDelta t}\frac{\left|\mathbf{r}^{\left(k\right)}\right|}{\left|\mathbf{\dot{r}}^{\left(k\right)}\right|}\mathbf{e}^{\left(k\right)}\cdot\mathbf{u}_{+}\right\} \left|\mathbf{r}_{+}^{\left(k+1\right)}\right|-\left\{ \frac{1}{4m\varDelta t}\frac{\left|\mathbf{r}^{\left(k\right)}\right|}{\left|\mathbf{\dot{r}}^{\left(k\right)}\right|}\mathbf{e}^{\left(k\right)}\cdot\mathbf{u}_{-}\right\} \left|\mathbf{r}_{-}^{\left(k+1\right)}\right|\nonumber \\
\qquad\qquad\qquad\qquad\qquad\qquad\qquad\qquad=\frac{\left|\mathbf{q}\right|}{m\left|\mathbf{r}^{\left(k\right)}\right|^{m-1}\sqrt{\left|\mathbf{\dot{r}}^{\left(k\right)}\right|}}+\left(1-\frac{1}{2m}\right)\left|\mathbf{r}^{\left(k\right)}\right|\,.\label{eq:general-reconstruct-equation}
\end{gather}
Cases of special interest are: a) $m=0$, $\mathbf{q}=\mathbf{u}\sqrt{\left|\dot{\mathbf{r}}\right|}$,
in this case the $\mathbb{L}^{2}$ metric expresses the length of
the contour%
\footnote{This is the case in the SRVF representation too.%
}, the reconstruction equations can be deduced from (\ref{eq:general-reconstruct-equation})
by multiplying both sides with $m$:
\begin{equation}
\left\{ \mathbf{e}^{\left(k\right)}\cdot\mathbf{u}_{+}\right\} \left|\mathbf{r}_{+}^{\left(k+1\right)}\right|-\left\{ \mathbf{e}^{\left(k\right)}\cdot\mathbf{u}_{-}\right\} \left|\mathbf{r}_{-}^{\left(k+1\right)}\right|=4\varDelta t\left(\left|\mathbf{q}\right|\sqrt{\left|\mathbf{\dot{r}}^{\left(k\right)}\right|}-\frac{1}{2}\left|\mathbf{\dot{r}}^{\left(k\right)}\right|\right)
\end{equation}
the coefficient matrix has special structure: the lack of diagonal
elements; b) $m=-\frac{1}{2}$, $\mathbf{q}=\mathbf{u}\sqrt{\frac{\left|\dot{\mathbf{r}}\right|}{\left|\mathbf{r}\right|}}$,
in this case the right hand side of (\ref{eq:linear-cost-solution-1})
is proportional to
\begin{equation}
1-\left|\mathbf{q}\right|\sqrt{\frac{\left|\mathbf{r}^{\left(k\right)}\right|}{\left|\dot{\mathbf{r}}^{\left(k\right)}\right|}},
\end{equation}
so at the solution this value becomes zero leading to homogeneous
equation system with the solution of identically zero $\left|\mathbf{r}\left(t\right)\right|$,
an obvious contradiction. The latter case is inherently singular as
already pointed out in \nameref{sec:Appendix-B}.

\bibliographystyle{plain}
\bibliography{refs}

\end{document}